\documentclass[psamsfonts,leqno,oneeqnum,onefignum]{siamltex}

\usepackage{amssymb}
\usepackage{amscd}
\usepackage{eucal}
\usepackage[dvips]{graphicx}

\newcommand{\beq}{\begin{equation}}
\newcommand{\eeq}{\end{equation}}

\newcommand{\ssum}[2]{\sum_{#1}^{#2} }
\newcommand{\psk}{\psi _k}
\newcommand{\intab}{\int_{a}^{b}}
\newcommand{\la}{\lambda}
\newcommand{\Tr}{\mbox{Tr\,}}
\newcommand{\N}{{\rm I\!N}}

\newcommand{\EB}{{E}}
\newcommand{\EA}{\mbox{E}}
\newcommand{\gb}{\overline{g}}
\newcommand{\SNR}{\mbox{SNR}}
\newcommand{\dB}{\mbox{dB}}
\newcommand{\binomial}[2]{(\!\!\begin{array}{c}
                                \raisebox{-.6ex}{$\scriptstyle #1$} \\
                                \raisebox{0.9 ex}{$\scriptstyle \!\! #2$}
                           \end{array}\!\! )}

\begin{document}

\title{On the Regularization of Fredholm Integral Equations of the First Kind\thanks{http://www.siam.org/journals/sima/29-4/30174.html}}

\author{
Enrico De Micheli\thanks{IBF - Consiglio Nazionale delle Ricerche,
Via De Marini 6, 16149 Genova, Italy ({\tt demicheli@ge.cnr.it}).}
\and
Nicodemo Magnoli\thanks{Dipartimento di Fisica -- Universit\`a di Genova,
Istituto Nazionale di Fisica Nucleare, sez. di Genova,
Via Dodecaneso 33, 16146 Genova, Italy
({\tt magnoli@ge.infn.it, viano@ge.infn.it}).}
\and
Giovanni Alberto Viano$^{\ddagger}$
}

\maketitle

\begin{abstract}
In this paper the problem of recovering a regularized
solution of the Fredholm integral equations of the first kind with
Hermitian and square-integrable kernels, and with
data corrupted by additive noise, is considered.
Instead of using a variational regularization of Tikhonov type, based on
a priori global bounds, we propose a method of truncation of eigenfunction
expansions that can be proved to converge asymptotically, in the sense of the
$L^2$--norm, in the limit of noise vanishing.
Here we extend the probabilistic counterpart of this procedure by
constructing a probabilistically regularized solution without assuming
any structure of order on the sequence of the Fourier coefficients
of the data. This probabilistic approach allows us to use
the statistical tools proper of time-series analysis, and in this way
we attain a new regularizing algorithm, which is illustrated by
some numerical examples. Finally, a comparison with solutions obtained by the means of the
variational regularization exhibits how some intrinsic limits of the
variational-based techniques can be overcome.
\end{abstract}



\pagestyle{myheadings}
\thispagestyle{plain}
\markboth{E. DE MICHELI, N. MAGNOLI, AND G. A. VIANO}{REGULARIZATION OF FREDHOLM INTEGRAL EQUATIONS}

\section{Introduction}
\label{introduction_section}
We consider the Fredholm integral equations of the first kind
\beq
\label{uno}
(Af)(x)=\intab K(x,y)f(y)\,dy = g(x)   ~~~~~~ (a \leq x \leq b)
\eeq
whose kernel $K(x,y)$ is supposed to be Hermitian and square integrable; i.e.,
\beq
\label{due}
K(x,y)=\overline{K(y,x)}
\eeq
and
\beq
\label{tre}
\intab \left \{ \intab | K(x,y) |^2 dx \right \}  dy < \infty.
\eeq
Then $A:L^2 (a,b) \rightarrow L^2 (a,b) $ is a self-adjoint compact operator.\\
For simplicity we shall suppose hereafter that the kernel $K$, the function $g$,
and the unknown function $f$ are real-valued functions; in addition, we assume
that the interval $[a,b]$ is a bounded and closed subset of the real line.

The Hilbert-Schmidt Theorem guarantees that the integral operator $A$
admits a set of eigenfunctions $\{ \psk \}_1^{\infty}$ and, accordingly,
a countably infinite set of eigenvalues $\{\la_k\}_{1}^{\infty}$. The
eigenfunctions form an orthonormal basis of the orthogonal complement
of the null space of the operator $A$ and therefore an orthonormal basis
of $L^{2}(a,b)$ when $A$ is injective. For the sake of simplicity only this
case will be considered, although this assumption can be easily relaxed with
slight technical modifications. The Hilbert-Schmidt theorem also guarantees
that $\lim_{k \rightarrow \infty} \lambda_k = 0$. Furthermore, we shall
suppose hereafter that the eigenvalues are ordered as follows:
$\la_1 > \la_2 > \la_3 > ....$ \\
In view of the Hilbert-Schmidt theorem we associate with the integral
equation (\ref{uno}) the following eigenfunction expansion:
\beq
\label{quattro}
f(x)=\ssum{k=1}{\infty} \left (\frac{g_k}{\lambda_k} \right ) \psi_k (x),
\eeq
where $g_k = (g,\psk )$, $((\cdot, \cdot)$ denoting the scalar product
in $L^2 (a,b))$. The series (\ref{quattro}) converges in the sense of $L^2$.

\emph{Remark.} If the support of the data does not coincide with that of the solutions,
i.e., $A:L^2(a,b) \rightarrow L^2(c,d)$ with $[a,b]$ different from $[c,d]$, the
problem can be worked out in terms of singular values and singular functions of the operator $A$
\cite{Bertero3}, and all of the following results can be easily reformulated.

In view of the fact that there always exists some inherent noise in the data,
instead of (\ref{uno}) we have to deal with the following equation:
\beq
\label{cinque}
Af+n=\overline{g}~~~~~~(\overline{g} = g + n),
\eeq
where $n$ represents the noise. Therefore, instead of expansion (\ref{quattro})
we have to consider the following expansion:
\beq
\label{sei}
\ssum{k=1}{\infty} \left (\frac{\overline{g}_k}{\lambda_k} \right ) \psi_k,
\eeq
where $\overline{g}_k = (\overline{g},\psk )$. Expansion (\ref{sei})
is generally diverging because $\overline{g} $ does not belong, in general,
to the range of the operator A. This is precisely a manifestation of
the ill-posed character of the Fredholm integral equation of the first kind.

Several methods of regularization have been proposed (see \cite{Engl,Groetsch,Hanke}
and references therein); all of them modify one of the elements of the triplet
$\{ A,X,Y \}$, where $A$ is the integral operator defined by (\ref{uno}),
whereas $X$ and $Y$ are, respectively, the solution and the data space (in our case
$X \equiv Y \equiv L^2 (a,b))$. Among these methods the procedure, which is
probably the most popular, consists in admitting only those solutions that belong
to a compact subset of the solution space $X$. In particular the famous method of Tikhonov
leads to the construction of ``regularizing operators'' by the minimization
of ``smoothing functionals''. In this latter functional the smoothing term is obtained
precisely by restricting the admitted solutions to a compact subset of the space $X$:
then the continuity of $A^{-1}$ follows from compactness. This restriction is
realized by the use of a priori bounds which can be written assuming some prior
knowledge of the solution. Therefore, in addition to the inequality
\beq
\label{sette}
\| Af - \overline{g} \| \leq \epsilon
\eeq
which corresponds to a bound on the noise ($\| \cdot \|$ denoting the norm in $L^2(a,b)$),
one also considers an a priori bound on the solution of the following form:
\beq
\label{otto}
\| Cf \|_{\cal Z} \leq \EB,
\eeq
where ${\cal Z}$ denotes the ``constraint space'' and, accordingly, $C$ is the
``constraint operator''.
From the bounds (\ref{sette}) and (\ref{otto}) we are led
to define the regularized solution as the minimum of the following functional:
\beq
\label{nove}
\Phi(f)=\left\| Af-\overline{g}\right\|^2 + \alpha^2 \left\| Cf\right\|_{\cal Z}^2,
~~~~~~~~~~\left ( \alpha = \left ( \frac{\epsilon}{\EB} \right ) \right ).
\eeq
In spite of several significant merits, this procedure is not free from defects.
Concerning the possibility of writing suitable a priori bounds on the solution,
we want to remark strongly that two different types of problems must be distinguished:
\begin{itemize}
\item[a)] synthesis problems;
\item[b)] inverse problems,
\end{itemize}
and to note that both are frequently solved by the use of Fredholm integral equations
of the first kind.
In the first class of problems, that basically consists in finding the source that produces
a prescribed effect (e.g. prescribed boundary values), the a priori bounds are intrinsic of
the problem itself, whereas this is not always the case for the second class.
As typical examples we can consider:
\begin{itemize}
\item[a$^\prime$)] the antenna synthesis;
\item[b$^\prime$)] the signal recovery.
\end{itemize}
The problem of the antenna synthesis consists in determining, within a certain degree of approximation,
the current intensity that generates a desired radiation pattern. It can be formulated in terms
of Fredholm equation of the first kind \cite{Maviano,Viano} and, consequently, it presents
the typical pathology of the ill-posed problems.
In this problem the a priori bound on the ohmic losses associated with the current intensity
is necessary and can be regarded as a natural constraint intrinsic of the
problem. Conversely, in the case of the signal recovery problem, the a priori bounds can be written only if
prior knowledge on the signal is assumed. Generally, it is possible to have some
a priori information regarding, for
instance, the support of the signal or requiring the function representing
the signal to be nonnegative. But even in these cases the prior knowledge
could be insufficiently specific to be peculiar of the function to be
reconstructed, and arbitrary, though reasonable, constraints must be added to solve the problem.
Strictly connected with this question there is the crux of the matter:
the practical choice of the regularization parameter $\alpha$ (see formula
(\ref{nove})) for a fixed $\overline{g}$, when the a priori bound (\ref{otto})
is unknown or it is not sufficiently precise. \\
Moreover, let us note that the functional (\ref{nove}) works as a filter
whose action is smoothing the Fourier components $\overline{g}_k$ for
high values of $k$. But it is easy to exhibit examples of signals
whose Fourier components are small, or even zero, for low values of $k$,
while the significant contributions of the signal are brought by those
components at intermediate values of $k$, which are smoothed out by the
action of the filter. In these situations the standard regularization
method fails,
showing that the only existence of the minumum of functional (\ref{nove})
does not guarantee the
bulk of the signal had been really recovered. This delicate point will be
illustrated with numerical examples in section \ref{numerical_analysis}.

We suggest a different approach which is based on the following observation:
for the moment, suppose that the moduli of the noiseless Fourier coefficients
$|g_k|$ are monotonically decreasing as $k$ increases; then, although
the formal series (\ref{sei}) diverges, nevertheless the effect of the error
remains limited in the beginning of the expansion, and there exists a point
(a certain value of $k$) where divergence sets in. Thus, the idea is to
stop the expansion at the point where it turns to diverge. This rough and
qualitative description can be put in rigorous form by proving that
even if the series (\ref{sei}) diverges, nevertheless it converges
(in the sense of $L^2$-norm) as $\epsilon$ (i.e., the bound on the noise)
tends to zero. This result, which has been proved by two of us (see \cite{Magnoli}),
does not give (except in very particular cases) a practical numerical
method for finding out the truncation point (i.e., the value of $k$) where
to stop expansion (\ref{sei}).
However, here we prove a probabilistic generalization of the results presented in \cite{Magnoli}
by removing the quite restrictive assumption that the Fourier coefficients
$|g_k|$ of the signal to be recovered are monotonically decreasing.
Compared to \cite{Magnoli} the significance of the new results is relevant.
First, the hypothesis made in \cite{Magnoli} on the order of the coefficients $|g_k|$
leads to a regularization procedure that essentially works as an ideal low-pass filter, and,
as previously discussed, this does not guarantee to recover correctly the signals whose
bulk is localized at intermediate frequencies.
Conversely, in this paper it will be shown how to construct a regularized solution
without assuming any kind of order on the coefficients $|g_k|$ by exploiting the tools supplied
by the information theory. This result will lead to a more effective
regularizing algorithm which is based on a suitable statistical analysis of the data and
whose main feature is indeed the frequency selectivity.
Second, from the application point of view, the hypothesis on the order
of the coefficients $|g_k|$ is too restrictive; thus, by removing it,
a much larger class of real signals can be practically analyzed.
These questions are precisely the contents of sections \ref{informatic_regularization}
and \ref{numerical_analysis}.
We will prove, indeed, in Section \ref{informatic_regularization} that it is possible
to split the noisy Fourier coefficients $\overline{g}_k$ into two classes:
\begin{itemize}
\item[i)] the Fourier coefficients $\overline{g}_k$ from which a significant
amount of information on $f_k = (f, \psi_k)$ can be extracted;
\item[ii)] the Fourier coefficients $\overline{g}_k$ that can be
regarded as random numbers because the noise prevails on the coefficients $g_k$.
\end{itemize}

In section \ref{numerical_analysis} it will be shown how it is possible to separate practically
the coefficients $\overline{g}_k$ into these two classes by the use of
statistical tools supplied by the so called ``time-series'' analysis.
Therefore, we can practically construct an approximation which converges
to the real solution, and furthermore we can have some confidence that the
bulk of the function $f$ has been effectively recovered.

The paper is organized as follows. In the first part of section \ref{variational_probabilistic_regularization}
a short sketch of the variational method based on the minimization of functional
(\ref{nove}) is given. This will be done in order to have explicitly the formulae
which will be used in section \ref{numerical_analysis}, where our procedure and the variational one
will be compared.
The second part of section \ref{variational_probabilistic_regularization}
is devoted to the probabilistic formulation of the regularization problem in a quite general setting.
In section \ref{informatic_regularization} we start illustrating the asymptotic convergence of the
eigenfunction expansion (in the sense of $L^2$-norm) as $\epsilon$ tends
to zero; then this result is reconsidered from the viewpoint of probability
and information theory. Here a key role will be played by the Bayes
formula: it will provide the various terms of our approximation, which will be
proved to be a probabilistically regularized solution of (\ref{uno}).
The first part of section \ref{numerical_analysis} is devoted to the discussion of the statistical tools
that are necessary for practically recovering the regularized solution from finite samples of noisy data. Finally,
some numerical examples are given in the second part of section \ref{numerical_analysis}.

\section{Variational and probabilistic regularization}
\label{variational_probabilistic_regularization}
\subsection{Variational regularization}
\label{variational_section}
After the classical book of Tikhonov and Arsenine \cite{Tikhonov} the literature
on the theory and applications of the variational regularization has been rapidly
growing (see, for instance, \cite{Groetsch}).
In order to compare our algorithm with this
classical one, some formulae and results of the variational regularization
will be here recalled here (see \cite{Bertero2,Bertero3,Miller1,Tikhonov} for proofs and details).

Let us characterize, first of all, the constraint operator $C$ and, accordingly, the constraint
space ${\cal Z}$. Let us take a constraint operator $C$ such that $C^\star C$ and
$A^\star A$ commute (this assumption does not restrict the theory and the applications
significantly \cite{Bertero2,Maviano}). Then, the space ${\cal Z}$ is composed
by those functions $f \in L^2(a,b)$ such that $\| Cf \|_{\cal Z}$ is finite; i.e.,
\beq
\label{dieci}
\| Cf \|_{\cal Z}^2 = (C^\star C f, f) =
\sum_{k=1}^\infty c_k^2 | f_k |^2 < \infty.
\eeq
Now we consider the ball
${\cal U_Z} = \{f \in {\cal Z} \, | \, \sum_{k=1}^\infty c_k^2 | f_k |^2 \leq \EB^2\}$,
and the restriction $A_0$ of the operator $A$ (see (\ref{uno})) to the ball
${\cal U_Z}$. Then, the following propositions can be proved.

\begin{proposition}
\label{pro:1}
If $\lim_{k\rightarrow \infty} c_k^2 = +\infty$ the operator $A_0^{-1}$ is continuous.
\end{proposition}

\begin{proposition}
\label{pro:2}
The functional $\Phi(f)$, with $\alpha= (\epsilon/\EB)$, has a unique minimum which is given by:
\beq
\label{dodici}
f_\star = [A^\star A + \left (\frac{\epsilon}{\EB}\right )^2 C^\star C]^{-1} A^\star \overline{g}.
\eeq
\end{proposition}

By expanding $\overline{g}$ in terms of $\psi_k$ (eigenfunctions of the operator $A$), we have:
\beq
\label{diciannove}
f_\star = \sum_{k=1}^\infty \frac{\lambda_k \overline{g}_k}{\lambda_k^2+c_k^2
\left (\frac{\epsilon}{\EB}\right )^2} \psi_k.
\eeq
Next, we have the following proposition.

\begin{proposition}
\label{pro:3}
The following limit holds true for any function $f$ satisfying the bounds $(\ref{sette})$ and $(\ref{otto})$:
\beq
\label{ventuno}
\lim_{\epsilon\rightarrow 0} \| f-f_\star \| = 0 ~~~~~ (E ~~\mbox{fixed}).
\eeq
\end{proposition}

In numerical computations it is often convenient to use truncated
approximations. For instance, one can derive from the smoothed solution (\ref{diciannove})
the following truncated approximation:
\beq
\label{ventisette}
f_\star^{(1)} = \sum_{k=1}^{k_\alpha} \frac{\overline{g}_k}{\lambda_k} \psi_k,
\eeq
where $k_\alpha$ is the largest integer such that
\beq
\label{ventotto}
\lambda_k \geq \left (\frac{\epsilon}{\EB}\right ) |c_k|.
\eeq

\begin{proposition}
\label{pro:4}
The following limit holds true for any function $f$ satisfying bounds $(\ref{sette})$ and $(\ref{otto})$:
\beq
\label{ventinove}
\lim_{\epsilon\rightarrow 0} \left\| f-f_\star^{(1)} \right\| = 0 ~~~~~ (E ~~\mbox{fixed}).
\eeq
\end{proposition}

In several problems a weaker a priori bound should be used by setting $C=I$ (the
identity operator). Therefore, instead of bound (\ref{otto}), we have
\beq
\label{trentasei}
\| f \| = \left ( \sum_{k=1}^\infty | f_k |^2 \right )^{1/2} \leq \EB.
\eeq
In this case the unique minimum of functional (\ref{nove}) is given by
\beq
\label{trentasette}
f_\star^{(2)} = \sum_{k=1}^\infty \frac{\lambda_k\overline{g}_k}{\lambda_k^2+
\left ( \frac{\epsilon}{\EB}\right )^2} \psi_k,
\eeq
and, accordingly, the following truncated approximation can be introduced:
\beq
\label{trentanove}
f_\star^{(3)} = \sum_{k=1}^{k_\beta} \frac{\overline{g}_k}{\lambda_k} \psi_k,
\eeq
where $k_\beta$ is the largest integer such that
\beq
\label{quaranta}
\lambda_k \geq \frac{\epsilon}{\EB}.
\eeq

Both $f_\star^{(2)}$ and $f_\star^{(3)}$ converge to $f$ as $\epsilon\rightarrow 0$ in a weak sense.
In fact, as shown in \cite{Miller1,Miller2}, the following proposition can be proved.

\begin{proposition}
\label{pro:5}
For any function $f$ which satisfies the bounds $(\ref{sette})$ and
$(\ref{trentasei})$, the following limits hold true:
\beq
\label{quarantanove}
\lim_{\epsilon\rightarrow 0} \left| \left (\left [f-f_\star^{(2)}\right ], v\right ) \right | = 0
~~~~~(\| v \| \leq 1, ~E~\mbox{fixed}),
\eeq
\beq
\label{cinquantuno}
\lim_{\epsilon\rightarrow 0} \left| \left (\left [f-f_\star^{(3)}\right ],v\right ) \right | = 0
~~~~(\| v \| \leq 1, ~E~\mbox{fixed}).
\eeq
\end{proposition}

\subsection{Probabilistic regularization}
\label{probabilistic_section}
Here we want to reconsider (\ref{cinque}) from a probabilistic point of view.
With this in mind we rewrite (\ref{cinque}) in the following form:
\beq
\label{sessantuno}
A\xi + \zeta =\eta,
\eeq
where $\xi$, $\zeta$ and $\eta$, which correspond to $f$, $n$ and
$\overline{g}$ respectively, are Gaussian weak random variables (w.r.v.)
in the Hilbert space $L^2 (a,b)$ \cite{Balakrishnan}.
A Gaussian w.r.v. is uniquely defined by its mean element and its
covariance operator; in the present case we denote by
$R_{\xi \xi}$, $R_{\zeta \zeta}$ and $R_{\eta \eta}$ the covariance operators
of $\xi$, $\zeta$ and $\eta$ respectively. Next, we make the following
assumptions:

\begin{itemize}
\item[I)] $\xi$ and $\zeta$ have zero mean; i.e. $m_\xi = m_\zeta = 0$;
\item[II)] $\xi$ and $\zeta$ are uncorrelated, i.e. $R_{\xi \zeta}$ = 0;
\item[III)] $R_{\zeta \zeta}^{-1}$ exists.
\end{itemize}

\noindent
The third assumption is the mathematical formulation of the fact that all
the components of the data function are affected by noise. As it is shown
by Franklin (see formula (3.11) of \cite{Franklin}), if the signal and the noise
satisfy assumptions I) and II), then
\beq
\label{sessantadue}
R_{\eta \eta}= A R_{\xi \xi} A^\star + R_{\zeta \zeta}
\eeq
and the cross-covariance operator is given by
\beq
\label{sessantatre}
R_{\xi \eta} = R_{\xi \xi}A^\star.
\eeq
We also assume that $R_{\zeta \zeta}$ will depend on a parameter $\epsilon$
that tends to zero when the noise vanishes; i.e.,
\beq
\label{sessantaquattro}
R_{\zeta \zeta} = \epsilon^{2} N,
\eeq
where $N$ is a given operator (e.g., $N=I$ for the white noise).

Now we are faced with the following problem.

\emph{Problem.} Given a value $\overline{g}$ of the w.r.v. $\eta$
find an estimate of the w.r.v. $\xi$.

A linear estimate of $\xi$ will be any w.r.v. $\xi_L=L\eta$, where $L:Y\rightarrow X$,
is an arbitrary linear continuous operator. Then, from a value $\overline{g}$ of $\eta$ one
obtains the linear estimate $L \overline{g}$ of the w.r.v. $\xi$. Now a measure of the reliability
of the estimator $L$ is given by
\beq
\label{sessantacinque}
\delta^2(\epsilon,v;L)=\EA \left\{ |(\xi-L\eta,v)|^2 \right\},~~(v \in X = L^2(a,b)),
\eeq
where $\EA\{\cdot\}$ denotes the expectation value. Then, we have the following proposition.

\begin{proposition}
\label{pro:8}
If the covariance operator $R_{\zeta\zeta}$ has a bounded inverse, then there
exists a unique operator $L_0$ that minimizes $\delta^2(\epsilon,v;L)$ for any $v\in X$,
and it is given by
\beq
\label{sessantasei}
L_0=R_{\xi \eta} R_{\eta \eta}^{-1}=R_{\xi \xi}A^\star \left [A R_{\xi \xi} A^\star + R_{\zeta \zeta}\right ]^{-1}.
\eeq
\end{proposition}

\begin{proof}
See \cite{Bertero1,Bertero2}.
\end{proof}

The w.r.v. $L_0\eta$ is called the best linear estimate of $\xi$, and, given a value
$\overline{g}$ of $\eta$, the best linear estimate $f_\star^{(4)}$ for the value of $\xi$ is
\beq
\label{sessantasette}
f_\star^{(4)}=\frac{R_{\xi \xi}A^\star}{A R_{\xi \xi} A^\star + R_{\zeta \zeta}} \overline{g},~~~~(A^\star=A).
\eeq
If $\xi$ and $L\eta$ have finite variance, then the global mean-square error may be defined as follows:
\beq
\label{sessantotto}
\delta^2(\epsilon,L)= \EA \left\{ \| \xi-L\eta\|^2\right\}.
\eeq
When the operator $L_0$ which minimizes (\ref{sessantacinque}) does exist, it also minimizes
the global error (\ref{sessantotto}) if $L_0\eta$ has finite variance; i.e., if
$\Tr (L_0 R_{\eta \eta} L_0^\star) < \infty$, then the following proposition can be proved.

\begin{proposition}
\label{pro:9}
If the following assumptions
\begin{itemize}
\item[{\rm i)}] $R_{\xi \xi}$ is an operator of trace class;
\item[{\rm ii)}] $R_{\zeta \zeta}=\epsilon^2 N$ has bounded inverse;
\item[{\rm iii)}] the equation $Af=0$, where $f\in\, \mbox{Range}\, \left (R_{\xi \xi}^{1/2}\right )$, has only the
trivial solution $f=0$
\end{itemize}
are satisfied, then the following limit holds true:
\beq
\label{sessantanove}
\lim_{\epsilon\rightarrow 0} \delta^2(\epsilon) = 0,
\eeq
where $\delta^2(\epsilon)=\inf_{L} \delta^2(\epsilon;L)$.
\end{proposition}

\begin{proof}
See \cite{Bertero1,Bertero2}.
\end{proof}

Let us note that $\delta^2(\epsilon)=\delta^2(\epsilon;L_0)$ when
$L_0$ does exist and is unique.

If we want to compare the probabilistic results obtained above with the variational ones,
which have been obtained by the use of eigenfunction expansions, we must expand $\xi$ and
$\zeta$ in terms of the eigenfunctions of the operator $A$ (i.e. $\{\psi_k\}_1^\infty$). Their
Fourier components are the random variables $\xi_k = (\xi, \psi_k)$ and $\zeta_k = (\zeta, \psi_k)$,
whose variances are given respectively by $\rho_k^2$ and $\epsilon^2 \nu_k^2$. Next, in addition
to the assumptions I)-III) made before, we make the following hypothesis in spite of the
fact that it turns out to be completely unrealistic (see section \ref{numerical_analysis}):
\begin{itemize}
\item[IV)] the Fourier components of $\xi$ are mutually uncorrelated as well as
the Fourier components of $\zeta$.
\end{itemize}
Therefore, if $R_{\zeta \zeta}^{-1}$ is bounded (i.e. $\sup_{k} (1/\epsilon^2 \nu_k^2)<\infty$),
then the operator $L_0$ exists and the best linear estimate (\ref{sessantasette}) can be written as
\beq
\label{settanta}
f_\star^{(4)}=\sum_{k=1}^\infty \frac{\lambda_k \rho_k^2}{\lambda_k^2 \rho_k^2+\epsilon^2 \nu_k^2} \overline{g}_k \psi_k.
\eeq
Finally, the quantities $\delta^2(\epsilon,v;L_0)$ and $\delta^2(\epsilon)$ become
\begin{eqnarray}
\label{settantuno}
\lefteqn{ \delta^2(\epsilon,v;L_0)=\EA \left\{| (\xi-L_0\eta,v)|^2\right\} =}
\hspace{1.3truecm} \nonumber \\[-1.5ex]
\\[-1.5ex]
& & = ([R_{\xi\xi}-L_0 R_{\eta\eta} L_0^\star]v,v)=
\epsilon^2\sum_{k=1}^\infty \frac{\rho_k^2 \nu_k^2}{\lambda_k^2 \rho_k^2+\epsilon^2 \nu_k^2} | v_k|^2
\hspace{0 truecm}\nonumber
\end{eqnarray}
and
\beq
\label{settantadue}
\delta^2(\epsilon)=\delta^2(\epsilon;L_0)=\Tr [R_{\xi\xi}-L_0 R_{\eta\eta} L_0^\star]=
\epsilon^2\sum_{k=1}^\infty \frac{\rho_k^2 \nu_k^2}{\lambda_k^2 \rho_k^2+\epsilon^2 \nu_k^2},
\eeq
and we have the following proposition.

\begin{proposition}
\label{pro:10}
The following statements hold true:
\begin{itemize}
\item[{\rm i)}] for any $v \in X$ ($X=L^2(a,b)$)
\end{itemize}
\beq
\label{settantatre}
\lim_{\epsilon\rightarrow 0} \delta^2(\epsilon,v;L_0) = 0,
\eeq
\begin{itemize}
\item[{\rm ii)}] if $\Tr R_{\xi\xi} < \infty$, then
\end{itemize}
\beq
\label{settantaquattro}
\lim_{\epsilon\rightarrow 0} \delta^2(\epsilon) = 0.
\eeq
\end{proposition}

\section{Information theory and regularization}
\label{informatic_regularization}
\subsection{Asymptotic convergence, in the $\mathbf{L^2}$-norm, of the eigenfunction expansion}
\label{asymptotic_section}
In the variational regularization, use is made of global a priori bounds (e.g., formulae
(\ref{otto}) or (\ref{trentasei})), which are the natural constraints in the case of synthesis
problems where the variational approach is certainly appropriate. But these global bounds
are not necessarily given in the case of inverse problems where the prior knowledge on the
solution can be, in several cases, rather poor. Moreover, in the truncated solutions derived by the
methods of variational regularization, the point at which to stop the expansion is obtained
by comparing the eigenvalues $\lambda_k$ with the ratio $(\epsilon/\EB)$ (i.e., formula (\ref{quaranta})),
or with $(\epsilon/\EB)|c_k|$ (see formula (\ref{ventotto})). In both cases
this approach appears quite unnatural from the viewpoint of the experimental or physical sciences,
whose methodology would rather suggest to stop the expansions
at the value $k_0$ of $k$ such that for $k > k_0$ the Fourier coefficients $g_k$ of
the noiseless data are smaller or at most of the same order of magnitude of $\epsilon$, and, consequently,
it is impossible to extract information from the corresponding coefficients $\overline{g}_k$.
With this in mind, and assuming that the noise is represented by
a bounded and integrable function $n(x)$ which satisfies the following condition:
\beq
\label{settantacinque}
\sup | n(x) | \leq \epsilon,~~~~x \in [a,b],
\eeq
the following results have been proved by two of us:

\begin{lemma}
\label{lem:2}
The following statements hold true:
\beq
\label{settantasei}
\sum_{k=1}^\infty \left (\frac{g_k}{\lambda_k}\right )^2 = \| f \|^2 = C_1 ~~~~~~ (C_1 = \mbox{constant}),
\eeq
\beq
\label{settantasette}
\sum_{k=1}^\infty \left ( \frac{\overline{g}_k}{\lambda_k}\right )^2 = + \infty
~~~~~\mbox{if}~ \overline{g} \not\in \mbox{Range}\,(A), ~~~~~~~~
\eeq
\beq
\label{settantotto}
\lim_{\epsilon\rightarrow 0} \overline{g}_k = g_k,~~\forall k. ~~~~~~~~~~~~~~~~~~~~~~~~~~~~~~~~~~~
\eeq

If $k_0(\epsilon)$ is defined by
\beq
\label{settantanove}
k_0(\epsilon) = \max \left\{m \in \N \, : \, \sum_{k=1}^m \left (
\frac{\overline{g}_k}{\lambda_k}\right )^2 \leq C_1\right\},
\eeq

then

\beq
\label{ottanta}
\lim_{\epsilon\rightarrow 0} k_0(\epsilon) = + \infty.
\eeq
\end{lemma}

\begin{proof}
See \cite{Magnoli}.
\end{proof}

Now we can introduce the following approximation
\beq
\label{ottantuno}
f_0^{(\epsilon)} = \sum_{k=1}^{k_0(\epsilon)}
\frac{\overline{g}_k}{\lambda_k} \psi_k
\eeq
and prove the following theorem.

\begin{theorem}
\label{the:1}
The following equality holds true:
\beq
\label{ottantadue}
\lim_{\epsilon\rightarrow 0} \left\| f-f_0^{(\epsilon)} \right\| = 0.
\eeq
\end{theorem}

\begin{proof}
See \cite{Magnoli}.
\end{proof}

If we consider a sequence of noisy data $\overline{g}$ which tends to $g$ for $\epsilon\rightarrow 0$ in the
sense of the $L^2$-norm (i.e., $\lim_{\epsilon\rightarrow 0} \| \overline{g} - g \| = 0$),
then $f_0^{(\epsilon)}$ will tend to $f$ as $\epsilon\rightarrow 0$ in the sense of the $L^2$-norm
(i.e., $\lim_{\epsilon\rightarrow 0} \| f_0^{(\epsilon)}-f\| = 0$).
In fact, since $\| \overline{g} - g \|^2=\sum_{k=1}^\infty | \overline{g}_k-g_k|^2$,
the $\lim_{\epsilon\rightarrow 0} \| \overline{g} - g \| = 0$ implies
that for any $k$, $\lim_{\epsilon\rightarrow 0} \overline{g}_k = g_k$, and
in view of Lemma \ref{lem:2} and Theorem \ref{the:1} it can be concluded that
$\lim_{\epsilon\rightarrow 0} \| f_0^{(\epsilon)}-f\| = 0$. Therefore,
from approximation (\ref{ottantuno}) we can derive an operator $\overline{B}$ defined by:
\beq
\label{ottantatre}
\overline{B} \overline{g} = \sum_{k=1}^{k_0(\epsilon)} \frac{\overline{g}_k}{\lambda_k} \psi_k,
\eeq
which continuously maps (i.e., preserving the convergence) the data $\overline{g}$ into the solution space $X$.
Thus, continuity has been restored without requiring compactness.

Two types of difficulties still remain:
\begin{itemize}
\item[a)] how to determine numerically the truncation point $k_0 (\epsilon)$, if
the norm of the function $f$ (i.e., the constant $C_1=\| f \|^2$) is unknown;
\item[b)] in any case the convergence of approximation (\ref{ottantuno}) is not sufficient to guarantee that the bulk
of the unknown function $f$ has been really recovered.
\end{itemize}
We can give a satisfactory answer to these questions only in very specific and peculiar cases, as we will
explain below. Suppose that the moduli of the Fourier coefficients $| g_k |$ are monotonically
decreasing for increasing values of $k$. Since $\overline{g}_k=g_k+n_k$, it turns out that at a
certain value $k_0$ of $k$ we have $| g_k | \simeq | n_k | \leq \epsilon$. The Fourier
coefficients of the noiseless data are of the same order of magnitude as the Fourier components of the
noise, and at this point we cannot extract any information from the noisy Fourier coefficients $\overline{g}_k$.
Let us now introduce the function $M(m)=\sum_{k=1}^m \left (\overline{g}_k/\lambda_k\right )^2$, whose
relevant properties are:
\begin{itemize}
\item[1)] It is an increasing function of $m$;
\item[2)] If $\epsilon$ is sufficiently small and the values of $| g_k |$ are monotonically
decreasing for increasing $k$, $M(m)$ presents a ``plateau'' when it reaches the value $C_1$. Indeed,
from formula (\ref{ottanta}) in Lemma \ref{lem:2} it
follows that $M(m)$ remains nearly constant when it attains the value $C_1$. An explicit numerical example
of this ``plateau'' is given in Figure \ref{figura_1}D in section \ref{numerical_analysis}.
\end{itemize}
This ``plateau'' corresponds to the order-disorder transition in the coefficients $\overline{g}_k$:
for $k < k_0(\epsilon)$ the data $g_k$ prevail on $n_k$ whereas for $k > k_0(\epsilon)$ the noise
components $n_k$ are larger or, at least, of the same order of magnitude of the noiseless data.
However, it must be remarked that in practical cases to single out the plateau which does really correspond
to the order-disorder transition in the coefficients $\overline{g}_k$ can be made difficult
by the presence of other spurious plateaux due to the erratic behavior of the noise.
Furthermore, if the coefficients $g_k$ are negligible for low values of $k$, and the
actual bulk of information is located only at intermediate values of $k$, there could be no numerical evidence of
such a plateau in spite of the fact that the convergence guaranteed by Theorem \ref{the:1} remains true.
Then we are forced to look for other methods that overcome these difficulties.
This issue will be investigated by means of probabilistic methods, as will be
illustrated in the next subsection.

\subsection{Bayes formula, information theory, and regularization}
\label{bayes_section}
Here our goal is to find a probabilistic extension of the result of Theorem \ref{the:1}
in which the assumption requiring the Fourier coefficients $|g_k|$ to be monotonically
decreasing will be removed.
In fact, we will show how to construct a regularizing solution from the noisy data,
disregarding the order of the coefficients $|g_k|$.
For this purpose, we turn (\ref{sessantuno}) into an infinite sequence of one-dimensional equations
by means of orthogonal projections:
\beq
\label{ottantaquattro}
\lambda_k \xi_k+\zeta_k=\eta_k,~~~(k=1,2,...),
\eeq
where $\xi_k = (\xi, \psi_k)$, $\zeta_k = (\zeta, \psi_k)$, $\eta_k = (\eta,\psi_k)$ are
Gaussian random variables. Here we retain assumptions I)-III) made in section \ref{probabilistic_section}, but
we remove assumption IV). In fact, there is no reason to assume that the basis
$\{\psi_k\}_1^\infty$ which diagonalizes the operator $A$ also diagonalizes the covariance
operators $R_{\xi \xi}$, $R_{\zeta \zeta}$, $R_{\eta \eta}$ \cite{Middleton}. Therefore, we can introduce
the variances $\rho_k^2 = (R_{\xi \xi} \psi_k, \psi_k)$,
$\epsilon^2 \nu_k^2 = (R_{\zeta \zeta} \psi_k, \psi_k)$,
$\lambda_k^2 \rho_k^2 + \epsilon^2 \nu_k^2 = (R_{\eta \eta} \psi_k, \psi_k)$,
without assuming that the Fourier components $\xi_k$ of $\xi$ (and analogously
also for $\zeta_k$ and $\eta_k$) are mutually uncorrelated. In view of the assumptions
I) and III) the following probability densities for $\xi_k$ and $\zeta_k$ can be assumed:
\beq
\label{ottantacinque}
p_{\xi_k} (x) = \frac{1}{\sqrt{2\pi}\, \rho_k } \exp \left \{ - \left ( \frac{x^2}{2\rho_k^2} \right ) \right\} ,
~~~(k=1,2, ...)
\eeq
and
\beq
\label{ottantasei}
p_{\zeta_k} (x) = \frac{1}{\sqrt{2\pi}\,\epsilon\nu_k} \exp\left\{ -
\left (\frac{x^2}{2\epsilon^2\nu_k^2}\right )\right\},
~~~(k=1,2, ...).
\eeq
By the use of the (\ref{ottantaquattro}) we can also introduce the conditional probability
density $p_{\eta_k} (y|x)$ of the random variable $\eta_k$ for fixed $\xi_k=x$, which reads
\begin{eqnarray}
\label{ottantasette}
p_{\eta_k} (y|x) &=& \frac{1}{\sqrt{2\pi}\,\epsilon\nu_k}
\exp\left\{-\frac{(y-\lambda_k x)^2}{2\epsilon^2\nu_k^2}\right\}\nonumber \\
&=& \frac{1}{\sqrt{2\pi}\,\epsilon\nu_k} \exp\left\{-\frac{\lambda_k^2}{2\epsilon^2\nu_k^2}
\left (x-\frac{y}{\lambda_k}\right )^2\right\}.
\end{eqnarray}
Now let us apply the Bayes formula that provides the conditional probability density of
$\xi_k$ given $\eta_k$ through the following expression:
\beq
\label{ottantotto}
p_{\xi_k} (x|y) = \frac{p_{\xi_k}(x) p_{\eta_k} (y|x)}{p_{\eta_k} (y)}.
\eeq
Thus, if a realization of the random variable $\eta_k$ is given by $\overline{g}_k$ (see
the formulation of the problem in section \ref{probabilistic_section}), formula (\ref{ottantotto}) becomes
\beq
\label{ottantanove}
p_{\xi_k} (x|\overline{g}_k) = A_k \exp\left\{-\frac{x^2}{2\rho_k^2}\right\}
\exp\left\{-\frac{\lambda_k^2}{2\epsilon^2\nu_k^2}\left (x-\frac{\overline{g}_k}{\lambda_k}\right )^2\right\}
~~~~(A_k = \mbox{const.}).
\eeq
Now the amount of information on the variable $\xi_k$ which is contained in the
variable $\eta_k$ can be evaluated. We have \cite{Gelfand}
\beq
\label{novanta}
J(\xi_k, \eta_k) = - \frac{1}{2} \log (1 - r_k^2),
\eeq
where
\beq
\label{novantuno}
r_k^2 = \frac{|\EA \left\{ \xi_k \eta_k\right\} |^2}{\EA\left\{ |\xi_k |^2 \right\} \EA\left\{ | \eta_k |^2\right\} } =
\frac{(\lambda_k \rho_k)^2}{(\lambda_k \rho_k )^2 + (\epsilon \nu_k )^2}.
\eeq
Thus,
\beq
\label{novantadue}
J (\xi_k, \eta_k) = \frac{1}{2} \log \left (1 + \frac{\lambda_k^2 \rho_k^2}{\epsilon^2 \nu_k^2} \right ).
\eeq
From equality (\ref{novantadue}) it follows that $J (\xi_k, \eta_k) < \frac{1}{2} \log 2$,
if $\lambda_k\rho_k < \epsilon \nu_k$. Thus, we are naturally led to introduce the following sets:
\beq
\label{novantatre}
{\cal I}_k = \left\{k\, :\, \lambda_k\rho_k \geq \epsilon\nu_k\right\},
\eeq
\beq
\label{novantaquattro}
{\cal N}_k = \left\{k\, :\, \lambda_k\rho_k < \epsilon\nu_k\right\}.
\eeq
Reverting to the conditional probability density (\ref{ottantanove}), it can be
regarded as the product of two Gaussian probability densities:
$p_1 (x) = A_k^{(1)} \exp\left\{-x^2/2\rho_k^2\right\}$
and
$p_2 (x) = A_k^{(2)} \exp\left\{-(\lambda_k^2/2\epsilon^2\nu_k^2)\left (x-(\overline{g}_k/\lambda_k)\right )^2\right\}$,
$(A_k=A_k^{(1)} \cdot A_k^{(2)})$, whose variances are respectively given by $\rho_k$
and $(\epsilon\nu_k/\lambda_k)$. Let us note
that if $k \in {\cal I}_k$, the variance associated with the density $p_2 (x)$ is smaller than the corresponding
variance of $p_1(x)$, and vice versa if $k \in {\cal N}_k$. Therefore, it is reasonable to consider as an acceptable
approximation of $\langle\xi_k\rangle$ the mean value given by the density $p_2(x)$ if $k \in {\cal I}_k$, or the mean
value given by the density $p_1(x)$ if $k \in {\cal N}_k$. We can write the following approximation:
\beq
\label{novantacinque}
\langle\xi_k\rangle = \left\{
\begin{array}{ll}
\frac{\strut\displaystyle \overline{g}_k}{\strut\displaystyle \lambda_k} &~~ (k \in {\cal I}_k), \\
0 &~~ (k \in {\cal N}_k).
\end{array}
\right .
\eeq
Consequently, given the value $\overline{g}$ of the w.r.v. $\eta$,
we are led to consider the following estimate of $\xi$:
\beq
\label{novantasei}
\widehat{B} \overline{g} = \sum_{k \in {\cal I}_k} \frac{\overline{g}_k}{\lambda _{k}} \psi_k.
\eeq
However, these are only heuristic considerations based on plausible arguments.
They will become rigorous statements only if it will proved that they lead to a solution $\widehat{B} \overline{g}$
which is probabilistically regularized.
For this purpose, the global mean-square error associated with the operator $\widehat{B}$, i.e.,
$\EA\left \{\|\xi - \widehat{B} \eta\|^2\right \}$, must be evaluated, and we have the following proposition.

\begin{proposition}
\label{pro:11}
\begin{itemize}
\item[{\rm i)}] If $\lim_{k\rightarrow\infty}(\lambda_k\rho_k/\nu_k) = 0$, then the set ${\cal I}_k$ is
finite for any fixed positive value of $\epsilon$;
\item[{\rm ii)}] assuming that the limit stated in {\rm i)} holds true, and, in addition, that $R_{\xi \xi}$ is an
operator of trace class, then the following relationship holds:
\end{itemize}
\beq
\label{novantasette}
\EA\left \{\|\xi - \widehat{B} \eta\|^2\right \} = \sum_{k\in {\cal N}_k} \rho_k^2 +
\sum_{k\in {\cal I}_k} \frac{\epsilon^2 \nu_k^2}{\lambda_k^2} < \infty.
\eeq
\end{proposition}

\begin{proof}
The proof of statement i) is obvious if we recall the definition of the
set ${\cal I}_k$ (formula (\ref{novantatre})). Statement ii) follows easily from the equality
\beq
\label{novantotto}
\EA\left\{\|\xi - \widehat{B} \eta\|^2\right \} = \Tr (R_{\xi\xi}-R_{\xi\xi} A^\star
\widehat{B}^\star - \widehat{B} A R_{\xi\xi} + \widehat{B} R_{\eta \eta} \widehat{B}^\star)
\eeq
and by the use of formulae (\ref{sessantadue}), (\ref{sessantaquattro}), and (\ref{novantasei}).
\end{proof}

In order to prove that approximation (\ref{novantasei}) is regularized, we need the following
auxiliary lemma.

\begin{lemma}
\label{lem:3}
Let $k_\gamma(\epsilon)$ be defined as follows:
\beq
\label{novantanove}
k_\gamma(\epsilon)= \max\,\left\{m\in \N \, :\, \sum_{k=1}^m
\left (\rho_k^2+\frac{\epsilon^2\nu_k^2}{\lambda_k^2}\right ) \leq \Gamma \right\},
\eeq
where $\Gamma = \Tr R_{\xi \xi}$ is finite. Then the following statements hold true:
\begin{eqnarray}
\label{cento}
\mbox{\rm i)} & ~~~ & \lim_{\epsilon\rightarrow 0} k_\gamma (\epsilon) = +\infty,\hspace{6truecm}\\
\label{centouno}
\mbox{\rm ii)} & ~~~ & \lim_{\epsilon\rightarrow 0} \left\{ \sum_{k=1}^{k_\gamma} \frac{\epsilon^2\nu_k^2}{\lambda_k^2} +
\sum_{k=k_\gamma+1}^\infty \rho_k^2 \right\} = 0.
\end{eqnarray}
\end{lemma}

\begin{proof}
i) Let $k_{\gamma_1}$ denote the sum $(k_\gamma+1)$. Then suppose that the limit
(\ref{cento}) does not hold. This latter assumption would imply that there exists a finite number $M$,
which does not depend on $\epsilon$, such that $k_{\gamma_1} < M$. Furthermore, this bound should remain true for any
sequence $\epsilon_i$ tending to zero. Then, we have the following inequalities:
\beq
\label{centodue}
\Gamma < \sum_{k=1}^{k_{\gamma_1}}\left (\rho_k^2+\frac{\epsilon^2\nu_k^2}{\lambda_k^2}\right ) \leq
\sum_{k=1}^{M} \left (\rho_k^2+\frac{\epsilon^2\nu_k^2}{\lambda_k^2}\right ).
\eeq
Now for any sequence $\epsilon_i$ tending to zero, we have
\beq
\label{centotre}
\Gamma < \sum_{k=1}^{M} \rho_k^2  \leq \sum_{k=1}^\infty \rho_k^2 = \Gamma,
\eeq
which is contradictory. Then, limit (\ref{cento}) holds.\\
ii) From $\sum_{k=1}^\infty \rho_k^2 = \Tr R_{\xi\xi}=\Gamma <\infty$, and
in view of statement i), it follows that
$\lim_{\epsilon\rightarrow 0} \sum_{k=k_{\gamma_1}}^\infty  \rho_k^2 =0$. Regarding the sum
$\sum_{k=1}^{k_\gamma}(\epsilon^2\nu_k^2/\lambda_k^2)$, we can proceed as follows.
From formula (\ref{novantanove}) we have
\beq
\label{centoquattro}
\sum_{k=1}^{k_\gamma}\frac{\epsilon^2\nu_k^2}{\lambda_k^2} +
\sum_{k=1}^{k_\gamma} \rho_k^2 \leq \Gamma.
\eeq
Then
\beq
\label{centocinque}
\sum_{k=1}^{k_\gamma} \frac{\epsilon^2\nu_k^2}{\lambda_k^2} \leq \Gamma - \sum_{k=1}^{k_\gamma} \rho_k^2 =
\sum_{k=k_{\gamma_1}}^\infty \rho_k^2,
\eeq
but in view of the fact that $\lim_{\epsilon\rightarrow 0} \sum_{k=k_{\gamma_1}}^\infty  \rho_k^2 =0$,
we have $\lim_{\epsilon\rightarrow 0}\sum_{k=1}^{k_\gamma} (\epsilon^2\nu_k^2/\lambda_k^2)=0$,
and statement ii) is proved.
\end{proof}

We can now prove the following theorem.

\begin{theorem}
\label{the:2}
If the covariance operator $R_{\xi\xi}$ is of trace class, and if the set
${\cal I}_k$ is finite (see Proposition $\ref{pro:11}$), then the following limit holds true:
\beq
\label{centosei}
\lim_{\epsilon\rightarrow 0} \delta^2(\epsilon,\widehat{B})=
\lim_{\epsilon\rightarrow 0} \EA\left \{ \|\xi - \widehat{B} \eta\|^2\right \} = 0;
\eeq
i.e., approximation $(\ref{novantasei})$ is probabilistically regularized.
\end{theorem}

\begin{proof}
In view of formula (\ref{novantasette}) in Proposition \ref{pro:11}, the proof
of equality (\ref{centosei}) reduces to the proof of the following limit:
\beq
\label{centosette}
\lim_{\epsilon\rightarrow 0} \left\{\sum_{k\in {\cal I}_k} \frac{\epsilon^2\nu_k^2}{\lambda_k^2} +
\sum_{k\in {\cal N}_k} \rho_k^2\right\} = 0.
\eeq
Regarding the first sum of (\ref{centosette}), we divide the set ${\cal I}_k$ into two subsets
defined by
\begin{eqnarray}
\label{centootto}
& {\cal I}_k^{(1)} = & \{k \in {\cal I}_k \, : \, k \leq k_\gamma\}, \\
\label{centonove}
& {\cal I}_k^{(2)} = & \{k \in {\cal I}_k \, : \, k > k_\gamma\};
~~~~~({\cal I}_k = {\cal I}_k^{(1)} \cup {\cal I}_k^{(2)}).
\end{eqnarray}
Accordingly, we can write
\beq
\label{centodieci}
\sum_{k\in {\cal I}_k}  \frac{\epsilon^2\nu_k^2}{\lambda_k^2} =
\sum_{k\in {\cal I}_k^{(1)}} \frac{\epsilon^2\nu_k^2}{\lambda_k^2}+\sum_{k\in {\cal I}_k^{(2)}}
\frac{\epsilon^2\nu_k^2}{\lambda_k^2}.
\eeq
Then $\sum_{k\in {\cal I}_k^{(1)}} (\epsilon^2\nu_k^2/\lambda_k^2) \leq \sum_{k=1}^{k_\gamma}
(\epsilon^2\nu_k^2/\lambda_k^2)$,
and in view of Lemma \ref{lem:3} (where we proved that
$\lim_{\epsilon\rightarrow 0}\sum_{k=1}^{k_\gamma} (\epsilon^2\nu_k^2/\lambda_k^2)=0$) it follows that
\beq
\label{centoundici}
\lim_{\epsilon\rightarrow 0} \sum_{k\in {\cal I}_k^{(1)}} \frac{\epsilon^2\nu_k^2}{\lambda_k^2} = 0.
\eeq
Regarding the term $\sum_{k\in {\cal I}_k^{(2)}} (\epsilon^2\nu_k^2/\lambda_k^2)$, since
$k \in {\cal I}_k$ then $\rho_k^2 \geq (\epsilon^2\nu_k^2/\lambda_k^2)$ and therefore
\beq
\label{centododici}
\sum_{k\in {\cal I}_k^{(2)}} \frac{\epsilon^2\nu_k^2}{\lambda_k^2} \leq
\sum_{k=k_{\gamma_1}}^\infty \rho_k^2.
\eeq
But, as we have seen in Lemma \ref{lem:3}, $\lim_{\epsilon\rightarrow 0} \sum_{k=k_{\gamma_1}}^\infty \rho_k^2 = 0$,
and consequently
\beq
\label{centotredici}
\lim_{\epsilon\rightarrow 0} \sum_{k\in {\cal I}_k^{(2)}} \frac{\epsilon^2\nu_k^2}{\lambda_k^2} = 0.
\eeq
We can conclude that $\lim_{\epsilon\rightarrow 0} \sum_{k\in {\cal I}_k} (\epsilon^2\nu_k^2/\lambda_k^2) = 0$.
Regarding the sum $\sum_{k\in {\cal N}_k}\rho_k^2$, we proceed in an analogous way by splitting the set
${\cal N}_k$ into two subsets defined by
\begin{eqnarray}
\label{centoquattordici}
& {\cal N}_k^{(1)} = & \{k \in {\cal N}_k \, : \, k \leq k_\gamma\}, \\
\label{centoquindici}
& {\cal N}_k^{(2)} = & \{k \in {\cal N}_k \, : \, k > k_\gamma\};
~~~~~({\cal N}_k = {\cal N}_k^{(1)} \cup {\cal N}_k^{(2)}).
\end{eqnarray}
Accordingly, we write
\beq
\label{centosedici}
\sum_{k\in {\cal N}_k}\rho_k^2 =
\sum_{k\in {\cal N}_k^{(1)}} \rho_k^2 + \sum_{k\in {\cal N}_k^{(2)}} \rho_k^2.
\eeq
If $k \in {\cal N}_k^{(1)}$ and by the use of inequality $\rho_k^2 < (\epsilon^2\nu_k^2/\lambda_k^2)$
(because $k \in {\cal N}_k$) we can write
\beq
\label{centodiciassette}
\sum_{k\in {\cal N}_k^{(1)}} \rho_k^2 \leq \sum_{k=1}^{k_\gamma} \frac{\epsilon^2\nu_k^2}{\lambda_k^2}.
\eeq
But in Lemma \ref{lem:3} we proved that
$\lim_{\epsilon\rightarrow 0} \sum_{k=1}^{k_\gamma} (\epsilon^2\nu_k^2/\lambda_k^2)=0$,
and therefore we have $\lim_{\epsilon\rightarrow 0} \sum_{k\in {\cal N}_k^{(1)}} \rho_k^2=0$.
Regarding the second term on the right-hand side of formula (\ref{centosedici}), we have
\beq
\label{centodiciotto}
\sum_{k\in {\cal N}_k^{(2)}} \rho_k^2 \leq \sum_{k=k_{\gamma_1}}^{\infty}\rho_k^2.
\eeq
But, again, $\lim_{\epsilon\rightarrow 0} \sum_{k=k_{\gamma_1}}^{\infty}\rho_k^2 = 0$, and then
$\lim_{\epsilon\rightarrow 0} \sum_{k\in {\cal N}_k^{(2)}} \rho_k^2 = 0$.
\end{proof}

\emph{Remarks.} i) It is worth it to notice that the proof of Theorem \ref{the:2} does not require
any type of order in the sum (\ref{novantasei}).
In fact, the only assumption that $\{\lambda_k\}$ is a strictly decreasing sequence does not
evidently imply that the terms $(\lambda_k\rho_k/\epsilon\nu_k)$ have any type of
monotonicity in $k$, and, consequently, the sum (\ref{novantasei}) cannot,
in general, be regarded as an ordered sum of terms up to a certain maximum value of $k$.
Thus, unlike the regularized solutions (\ref{diciannove}), (\ref{ventisette}),
(\ref{trentasette}), (\ref{trentanove}) and also (\ref{ottantatre}), $\widehat{B}\gb$
features frequency selectivity, which is obtained by evaluating the information
content of the noisy Fourier coefficients. \\
ii) Notice that the estimate (\ref{novantasei}) associated with the operator $\widehat{B}$ represents
a probabilistically regularized solution, in the sense of the formula (\ref{centosei}),
even if, in general, it does not minimize the global mean-square error (\ref{sessantotto}).

~

At this point in order to apply the results of this section, statistical methods that allow for
splitting the coefficients
$\overline{g}_k$ into the sets ${\cal I}_k$ and ${\cal N}_k$ must be investigated. These methods will be illustrated
in the next section.

\section{Numerical analysis: the regularizing algorithm}
\label{numerical_analysis}
\subsection{The correlation function of the noisy data}
\label{correlation_section}
The application of the results of the previous section to a Fredholm equation
of the first kind would involve using statistical tools for the determination of the two sets
${\cal I}_k$ and ${\cal N}_k$.
In this section this issue is discussed and the basic steps of a numerical algorithm for constructing
the regularized solution $\widehat{B}\gb$ from the noisy data $\gb$  are outlined.
For simplicity we shall work throughout only with data corrupted by white noise.
However, provided the independence assumption between $\xi$ and $\zeta$,
more general cases involving ``colored'' noise
could be treated by using suitable methods, for instance, ``prewhitening'' transformations \cite{Box},
whose discussion is beyond the scope of this section.
Here our goal is to show that statistical estimates of the amount of information carried by
the Fourier coefficients $\gb_k$ can be sufficient to construct a satisfactory regularized solution.
Furthermore, the direct comparison of the numerical results clearly evidentiates how some inherent limitations
of the variational regularization scheme are overcome.

Following the analysis of the previous section, we are now faced with the problem of separating
the Fourier coefficients $\gb_k$ into two classes; one containing all the Fourier coefficients
of the noisy data which are correlated, the other containing the $\gb_k$ that can be regarded
as random numbers.
This task can be achieved by computing the correlation function of the random variables
$\eta_k$: i.e., the probabilistic counterpart of the coefficients $\gb_k$:
\begin{equation}
\label{IV-1}
\Delta_{\eta}(k_1, k_2) = \frac{E\{ [\eta_{k_1}-E\{\eta_{k_1}\}]  [\eta_{k_2}-E\{\eta_{k_2}\}]\}}
{E\{[\eta_{k_1}-E\{\eta_{k_1}\}]^2\}^{1/2} E\{[\eta_{k_2}-E\{\eta_{k_2}\}]^2\}^{1/2}}.
\end{equation}
In practice, just a finite realization $\{\gb_k\}_1^N$ of the random variables $\eta_k$ is available, from
which estimates $\delta_{\gb}$ of the autocorrelations can be obtained by regarding the data
$\{\gb_k\}_1^N$ as a finite length
record of a stationary random normal series. In principle the assumption of stationarity of the series
$\{\eta_k\}$ is not correct because in general the moments of the random variables $\eta_k$ will depend on $k$,
but from the practical point of view this is usually the only possible chance. In fact,
in many areas of application, it is difficult or even impossible to have multiple independent realizations
$\{\gb_k\}_1^N$ of the process $\{\eta_k\}$, so
estimates of ensemble averages cannot be computed. Thus, we are forced to introduce the working hypothesis
that the process $\{\eta_k\}$ is stationary in wide sense \cite{Doob}, that is
$\Delta_{\eta}(k_1, k_2) = \Delta_{\eta}(k_1 - k_2)$,
and compute the estimates of the autocorrelation coefficients
by means of the ergodic relation between ensemble and {\it time} (i.e., the index $k$ in our case) averages.
Of course, such a restriction can be removed whenever many independent sets of data $\{\gb_k\}_1^N$ would be
available for evaluating ensemble averages. Anyway, we will see later in the discussion of the algorithm how
an ambiguity in the reconstruction of the regularized solution $\widehat{B}\gb$
due to the assumed invariance for $k$-translation of $\{\eta_k\}$ will be removed.

\noindent
A number of estimators of the autocorrelation function have been suggested by statisticians
and their properties are discussed in detail in \cite{Jenkins}.
An estimate which is widely used by statisticians, and in the following examples as well, is given by
\begin{equation}
\label{IV-2}
~~~~\delta_{\gb} (n) =
\frac{\strut\displaystyle \sum_{k=1}^{N-n} (\gb_k - \langle\gb_k\rangle) (\gb_{k+n} - \langle\gb_{k+n}\rangle)}
     {\left\{\strut\displaystyle \sum_{k=1}^{N-n} (\gb_k - \langle\gb_k\rangle)^2 \sum_{k=1}^{N-n}
      (\gb_{k+n} - \langle\gb_{k+n}\rangle)^2\right\}^{1/2}},~~~~
n = 0, ..., N-1,
\end{equation}
where
\begin{equation}
\label{IV-3}
\langle\gb_k\rangle = \frac{1}{N-n} \sum_{k=1}^{N-n} \gb_k;~~~~~~ \langle\gb_{k+n}\rangle = \frac{1}{N-n} \sum_{k=1}^{N-n} \gb_{k+n}.
\end{equation}
\noindent
Equation (\ref{IV-2}), which is based on the scatter diagram of $\gb_{k+n}$ against $\gb_k$ for
$k = 1, .., N-n$, represents the maximum likelihood estimate of the autocorrelation coefficients
of two random variables $\eta_k$ and $\eta_{k+n}$ whose joint probability distribution function
is bivariate normal.

In order to identify the structure of the series $\{\bar{g}_k\}_1^N$ so that we can separate correlated components
from the random ones, it is necessary to have a crude test
on whether $\delta_{\gb} (n)$ is effectively zero.
It has been shown by Anderson \cite{Anderson} that the distribution of an estimated
autocorrelation coefficient, whose theoretical value is zero, is approximately normal.
Thus, on the hypothesis that the theoretical autocorrelation $\Delta_{\eta}(n) = 0$, the
estimate $\delta_{\gb}(n)$ divided by its standard error $\sigma_\delta(n)$ will be
approximately distributed as a unit normal deviate. This fact may be used to provide a rough guide as to
whether theoretical autocorrelations are essentially zero.
To this purpose it is usually sufficient to remember that, for normal distribution,
deviations exceeding two standard errors in either direction have a probability of
about $0.05$, so that the 95\% confidence interval of the estimate is approximately
$\delta_{\gb}(n) \pm 1.96 \,\sigma_\delta(n)$.

Estimated autocorrelations can have rather large variances and can be highly correlated
with each other \cite{Bartlett,Fuller}, so that care is required in the interpretation of individual autocorrelations.
In particular, moderately large estimated autocorrelations can occur after the theoretical
autocorrelation function has damped out and, in any case, it must be considered
that an estimated autocorrelation function always exhibits less damping than the theoretical one,
as the estimated autocorrelations are inflated by sampling fluctuations (see also the following Example 1).
Thus, in order to avoid a purely empirical analysis of the autocorrelations, it is necessary to assume a
rough model of the series that allows to evaluate the order of magnitude of the sampling errors $\sigma_\delta(n)$
associated to the autocorrelation estimator.

According to the discussion of Section \ref{bayes_section}, since we are expected to find
the set ${\cal I}_k$ to be finite, we are
also expected that the autocorrelation function $\Delta_{\eta}(n)$ will vanish beyond a certain lag $n_0$.
Thus, in what follows, it will be assumed that there exists an index $n_0$ such that $\Delta_{\eta}(n)=0$ for $n>n_0$.
In this case, if the record length $N$ is large enough (i.e., such that $O(1/N^2)$ terms
can be neglected), use can be made of the Bartlett's approximate expression for
the variance of the estimated autocorrelations of a stationary normal process \cite{Bartlett}:
\begin{equation}
\label{IV-4}
{\rm var}\, [\delta_{\gb}(n)] \sim \frac{1}{N-n} \left\{1+2\sum_{v=1}^{n_0}\Delta_{\eta}^2(v)\right\},
~~~~\mbox{for}~~ n > n_0.
\end{equation}
To use (\ref{IV-4}) in practice, the estimated autocorrelations $\delta_{\gb}$ are substituted for
the theoretical ones $\Delta_{\eta}$, and in this case we shall refer to the square root of (\ref{IV-4})
as the {\em large-lag} standard error $\sigma_\delta (n;n_0)$ \cite{Box}.

The index $n_0$ is actually recovered in a recursive way through an hypothesis generation-verification procedure.
Starting from the assumption that the series is completely random, i.e., $n_0=0$, the standard error
$\sigma_\delta (n;0)$ is computed and the first index $\overline{n}>0$ such that
$|\delta_{\gb}(\overline{n})| > 1.96\,\sigma_\delta (n;0)$ is searched for. If there exists such an
index $\overline{n}$, it
becomes the new candidate to be $n_0$, i.e., we set $n_0=\overline{n}$, $\sigma_\delta (n;n_0)$ is computed,
and again it is tested whether the series is compatible with the hypothesis that $\Delta_{\eta}(n)=0$ for $n>n_0$.
The whole procedure is repeated until no new index $\overline{n}$ is found. Formally, $n_0$ is then
defined as
\begin{equation}
\label{IV-5}
n_0 = \max \, \{\overline{n} \geq 0 \, : \, \forall n \in (\overline{n},N-1], \,
\mid \delta_{\gb}(n)\mid < 1.96~ \sigma_\delta (n,\overline{n})\}.
\end{equation}
The set ${\bf Q}$ of the lags corresponding to autocorrelation values that are effectively different from zero
and, consequently, indicating lack of randomness of the coefficients $\gb_k$, is defined as:
\begin{equation}
\label{IV-6}
{\bf Q} = \{0 < n \leq n_0\, : \, \mid \delta_{\gb} (n) \mid > 1.96~\sigma_\delta(n,0)\}.
\end{equation}
Let $N_c$ be the number of elements of $\bf Q$.

As previously discussed, as a consequence of the inevitable assumption of stationarity of the process $\{\eta_k\}$,
the Fourier coefficients $\gb_k$ that are correlated cannot be determined in a unique way from the set $\bf Q$.
In fact, an integer $n_i \in {\bf Q}$ just indicates a strong correlation between at least two Fourier coefficients
$n_i$ apart. This means that, in principle, any couple $(\gb_{k_i}, \gb_{k_i+n_i})$ for any integer
$1 \leq k_i \leq (N-n_i)$
could have generated such a strong correlation at the lag $n_i$.
Thus, from the set $\bf Q$ we can construct $N_c$ families $F_i$ defined as
\begin{equation}
\label{IV-7}
F_i = \left\{(\gb_{k_i}, \gb_{k_i+n_i})\right\}_{k_i=1}^{(N-n_i)},~~~i=1,\, ...,\, N_c
\end{equation}
from which the couples of coefficients $\gb_k$ that are likely to be correlated can be selected.
In theory, that is for $N\rightarrow\infty$, the $N_c$ indices $k_i$ and the $N_c$ elements $n_i \in {\bf Q}$ are
mutually dependent. In fact, any two coefficients $\gb_{k_\alpha}, \gb_{k_\beta}$ which are
selected from the families $F_i$
must satisfy the pairwise compatibility conditions requiring $|k_\alpha-k_\beta|\in {\bf Q}$.
Or, in other words, it can be seen that, given the set ${\bf Q}$, the number $N_{\cal I}$ of
admissible Fourier coefficients
$\gb_k$ is combinatorially constrained to be
\begin{equation}
\label{IV-8}
\frac{1}{2}\,(1+\sqrt{1+8N_c}) \leq N_{\cal I} \leq N_c+1.
\end{equation}
The left inequality in (\ref{IV-8}) follows directly from the observation that the maximum number of correlations
among $N_{\cal I}$ coefficients is $\binomial{N_{\cal I}}{2}$, then $N_c \leq \binomial{N_{\cal I}}{2}$,
whereas the right inequality
expresses that at least $(N_{\cal I}-1)$ distinct correlations can be computed among $N_{\cal I}$
coefficients (i.e. $N_c \geq N_{\cal I}-1$).
For instance, if $N_c = 2$, we have from inequalities (\ref{IV-8}) that there need to be $N_{\cal I}=3$
coefficients $\gb_k$ to construct
the set ${\bf Q}$, or, referring to (\ref{IV-7}), that the two indices $k_1$ and $k_2$ must coincide, i.e.,
$k_1\equiv k_2 \geq 1$.
In any case, the compatibility conditions are not sufficient to constraint in a unique way the
selection of the coefficients $\gb_k$ and, consequently, the construction of the regularized
solution.\\
In practice, that is when the record length $N$ is finite and particularly when the
signal-to-noise ratio (SNR) of the data $\gb$ is small, the compatibility constraints
cannot be assumed to be satisfied. In fact, because of the sampling fluctuations in the
estimates $\delta_{\gb} (n)$, some correlations which are actually different from zero could be incorrectly
detected by the procedure discussed above. However, we shall see later in the discussion
of the numerical examples how the compatibility constraints can provide us with a confidence
check on the reliability of the regularized solution $\widehat{B}\gb$.

In order to recover in a unique way from the set ${\bf Q}$ the Fourier coefficients that are likely
to be correlated, we adopt the following criterion
suggested by the definition itself of the autocorrelation function: for any $n_i \in {\bf Q},\, i=1,...,N_c$,
we select the pair $(\gb_{k^\star_i}, \gb_{k^\star _i+n_i})$ giving the maximum contribution to
the autocorrelation estimate
$\delta_{\gb} (n_i)$; i.e., we define $k^\star_i$ as
\begin{equation}
\label{IV-9}
k^\star_i = \arg \max_{k\in [1,N-n_i]}\,\{|\gb_k\,\gb_{k+n_i}|\},~~~i=1,\, ...,N_c,
\end{equation}
and, accordingly, we can define the set of frequencies ${\cal I}_k$ exhibiting correlated Fourier coefficients as
\begin{equation}
\label{IV-10}
{\cal I}_k = \{k^\star _i\}_1^{N_c} \cup \{k^\star_i+n_i\}_1^{N_c},
\end{equation}
where each element of ${\cal I}_k$ is counted only once.

\subsection{Numerical examples}
\label{example_section}
Throughout this section we shall consider as a sample problem the integral equation (\ref{uno}) with kernel
\begin{equation}
\label{IV-11}
K(x,y) = \left \{
\begin{array}{ll}
(1-x)\,y & ~~~\mbox{if ~~ $0 \leq y \leq x \leq 1$}, \\
x\,(1-y) & ~~~\mbox{if ~~ $0 \leq x \leq y \leq 1$}
\end{array}
\right.
\end{equation}
whose eigenfunctions and eigenvalues are, respectively,
\begin{eqnarray}
\label{IV-12}
\psi_k (x) & = & \sqrt 2 \,\sin (k \pi x), \\
\lambda_k  & = & \frac{1}{k^2 \pi^2}.
\end{eqnarray}

\noindent
The data $g(x)$ have been noised by adding white noise $n(x)$,
simulated by computer generated random numbers uniformly distributed in the interval $[-\epsilon,\epsilon]$
(see also \cite{Scalas} for a very preliminary numerical analysis of this problem).
The examples shown hereafter differ for the choice of the input signal $f(x)$ and for the values of the noise
boundary $\epsilon$,
whereas the performances of the algorithm are evaluated by direct comparison of the reconstructed signal
with the true signal $f(x)$.
In every example reported here, the approximations obtained through the variational scheme (see section
\ref{variational_section}) are computed by setting
the constraint operator $C$ such that $c_k=k,~(k=1,2,...)$, the parameter $\epsilon$ corresponding to the
boundary on the noise equal to the dispersion of the noise $D_\epsilon$ (see (\ref{sette})),
and the boundary $E$ on the solution equal to the norm of the unknown function, i.e., $E=\|f(x)\|$ (see (\ref{otto})).

\begin{figure}[tbh]
\begin{center}
\includegraphics[width=9cm]{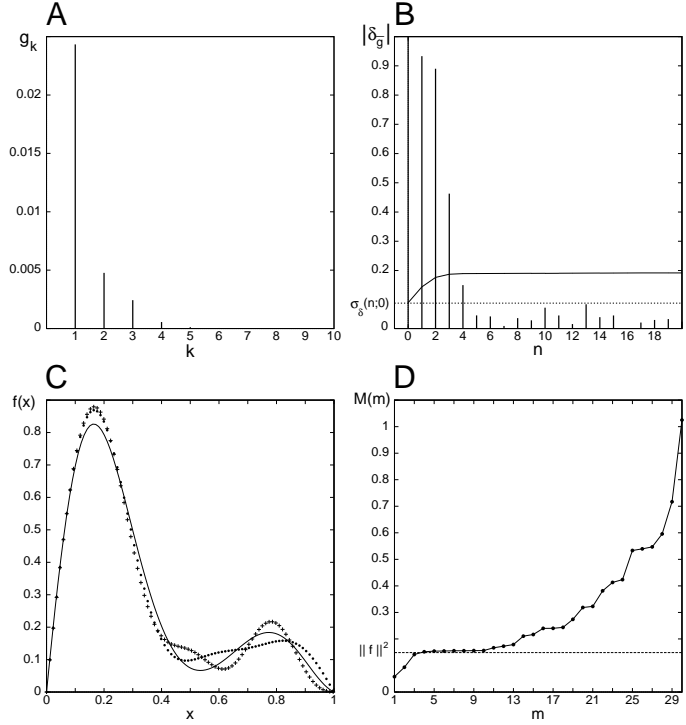}
\caption{\label{figura_1} Example 1: $f_1(x)=(1-x)\sin(3\sin(3x))$, $\epsilon = 10^{-4}$, $\SNR \simeq 25.7\dB$, $N = 512$.
(A) Noiseless Fourier coefficients $g_k$.
(B) Modulus of the autocorrelation function. The horizontal dotted straight line indicates the 95\% confidence limit
$1.96\,\sigma_\delta(n;0)$ for a purely random sequence. The solid curved line indicates the confidence limit
$1.96\,\sigma_\delta(n;3),\,n>3$. From the analysis of $\delta_{\overline{g}}(n)$ we have ${\bf Q} = \{1,2,3\}$
and ${\cal I}_k = \{1,2,3,4\}$.
(C) Regularized solutions. The solid line represents the actual solution $f_1(x)$. The dots represent the
reconstruction $\widehat{B}\overline{g}$. The crosses represent the variational solution $f_{\star}^{(1)}$
obtained by using $c_k = k$; $k_\alpha = 8$ (see equations (\protect\ref{ventisette}) and (\protect\ref{ventotto})).
(D) Plot of the function $M(m) = \sum_{k=1}^m\, (\overline{g}_k/\lambda_k)^2$. Notice that the value of $M(m)$
corresponding to its first plateau, i.e., approximately for $4 \leq m \leq 10$, is about the squared
norm of the true solution.
}
\end{center}
\end{figure}

In Figure \ref{figura_1}, the analysis of the sample function
$f_1(x)=(1-x)\sin(3\sin(3x))$ with noise boundary $\epsilon=10^{-4}$ is summarized.
The global $\SNR$,
defined as the ratio of the mean power of the noiseless data to the noise variance, was $\SNR \simeq 25.7 \dB$.
The function $f_1(x)$ is characterized by having the bulk of information localized in the first few values of $k$
(see the related noiseless coefficients $g_k$ in Figure \ref{figura_1}A) so that we expect that also a variational
solution could provide a satisfactory reconstruction of the input signal.
Figure \ref{figura_1}B shows the behavior of the autocorrelation function $\delta_{\gb}(n)$ along with
the two lines indicating the statistical confidence limits we used to discriminate whether the
autocorrelations are essentially null.
The dashed horizontal straight line represents the threshold that we would have under the
hypothesis of purely random sequence $\{\gb_k\}$, whereas the solid line represents
the threshold corresponding to the model of autocorrelation function of ideal damped type. In this example we found
$n_0 = 3$, ${\bf Q}=\{1,2,3\}$, and the autocorrelation at $n=4$ was rejected in spite of its quite large value
(see formula (\ref{IV-6})).
The direct inspection of the values of $\epsilon\nu_k$ and $g_k$ in repeated realizations showed that for $k=5$ the noise
was usually larger than the Fourier coefficient, confirming hence the result that the autocorrelation
$\delta_{\gb}(4)$ was
abnormally inflated by the large autocorrelations at $n=1,2,3$. According to the criteria (\ref{IV-9}) and (\ref{IV-10}),
the set of frequencies whose corresponding Fourier coefficients exhibit strong correlations is ${\cal I}_k=\{1,2,3,4\}$.
It is worth noticing that in this case the elements of ${\cal I}_k$ satisfy all the compatibility constraints,
i.e. any difference between elements of ${\cal I}_k$ belongs to ${\bf Q}$, and $N_{\cal I}$ satisfies
constraints (\ref{IV-8}).
This complete cross-consistency between ${\bf Q}$ and ${\cal I}_k$ gives a high level of confidence in the
result of the whole analysis.
In Figure \ref{figura_1}C the true function $f_1$ (solid line), the regularized solution
$\widehat{B}\gb$ (crosses) and the regularized function $f_{\star}^{(1)}$
(dots) are compared.
The truncation point of $f_{\star}^{(1)}$, obtained through the criterion (\ref{ventotto}), was $\alpha=8$.
Figure \ref{figura_1}C shows how in this case both regularization methods lead to comparable results, which
are quite satisfactory approximations of the ``unknown'' function $f_1$.
The plot of the function $M(m)$, displayed in Figure \ref{figura_1}D confirms the correctness of the
two approximations. In fact, it clearly exhibits a ``plateau'', ranging from about $m=3$ to $m=10$,
that corresponds to the order-disorder transition of the coefficients $\gb_k$. Then it could be
argued that for any truncation point belonging to this ``plateau'' the truncated approximation will
hold coefficients $\gb_k$ whose information content is not completely obscured by the noise.
In every example discussed here, the regularized solutions $f_\star^{(2)}$ and $f_\star^{(3)}$ (see
(\ref{trentasette}) and (\ref{trentanove})) have also been considered, providing
in all cases worse results (not plotted).
\begin{figure}[tb]
\begin{center}
\includegraphics[width=8cm]{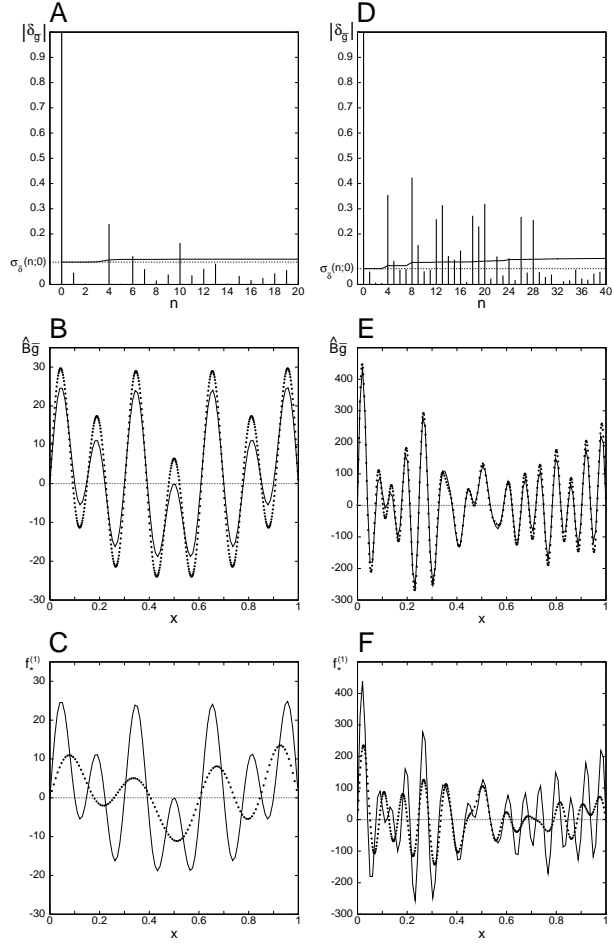}
\caption{\label{figura_2}
Example 2: $f_2(x) = 5\sin(3\pi x)+10\sin(7\pi x)+15\sin(13\pi x)$, $\epsilon = 3\,10^{-3}$,
$\SNR \simeq 0.54 \dB$, $N = 512$.
(A) Modulus of the autocorrelation function. ${\bf Q} = \{4,6,10\}$, ${\cal I}_k = \{3,7,13\}$.
(B) Comparison between the actual solution $f_2(x)$ (solid line) and the regularized solution
$\widehat{B}\overline{g} (x)$ (dots).
(C) Comparison between the actual solution $f_2(x)$ (solid line) and the approximated solution
$f_{\star}^{(1)}(x)$ with $k_\alpha = 9$ (see criterion (\protect\ref{ventotto})).
(D) Example 3: Modulus of the autocorrelation function. $f_3(x) = \sum_{j=1}^{10} \, a_j\,\sin(k_j \pi x)$,
with $a_j = \{17,23,27,33,43,55,68,70,77,81\}$ and $k_j=\{5,9,13,17,18,23,24,25,31,33\}$.
$\epsilon = 10^{-3}$, $\SNR \simeq 9.79 \dB$, $N = 1024$;
${\bf Q} = \{4,5,8,9,12,13,14,15,16,18,19,20,22,24,26,28\}$,
${\cal I}_k = \{5,9,13,17,18,23,24,25,31,33\}$.
(E) Comparison between the actual solution $f_3(x)$ (solid line) and the regularized solution
$\widehat{B}\overline{g}(x)$ (dots).
(F) Comparison between the actual solution $f_2(x)$ (solid line) and the approximated solution
$f_{\star}^{(1)}(x)$ with $k_\alpha=27$.
}
\end{center}
\end{figure}

The second and third examples, shown in Figure \ref{figura_2}, are quite simple but a little tricky, and show
the deep differences between our approach and the variational one. They consist of a finite linear
combination of, respectively,
3 and 10 basis functions $\psi_k$ (see the legend for numerical details), and, indeed, they have been chosen as
typical signals in which the bulk of the information is not grouped in a single block of consecutive low frequencies.
In these cases, setting global constraints on the solution, such as in the variational methods, leads
inevitably to a failure,
which is clearly evident from Figure \ref{figura_2}C,F, since the lack of selectivity necessarily causes the
regularized solution $f_{\star}^{(1)}$ to contain pure noisy components. On the contrary, the selectivity
achieved through the analysis of the autocorrelation function overcomes this limit.
In both examples the analysis of the autocorrelation function (see Figure \ref{figura_2}A,D) led to the correct
selection of the components that carry information in spite of the quite small $\SNR$
(in the Example 2, $\SNR \simeq 0.55 \dB$).
Referring to the Example 2 depicted in Figure \ref{figura_2}A,B,C, it can be observed that all the compatibility
constraints are indeed satisfied; however, it is worth to remark that, because of the sampling fluctuation of the
estimates $\delta_{\gb}(n)$, the autocorrelation $\delta_{\gb}(6)$ was not always detected in different realizations
of the noisy data $\{\gb_k\}_1^N$.
In these cases the set ${\cal I}_k$, computed from the set ${\bf Q}=\{4,10\}$ missing $n=6$, is still correct, i.e.,
${\cal I}_k=\{3,7,13\}$, even though one compatibility constraint is not fulfilled.

\begin{figure}[tbh]
\begin{center}
\includegraphics[width=9cm]{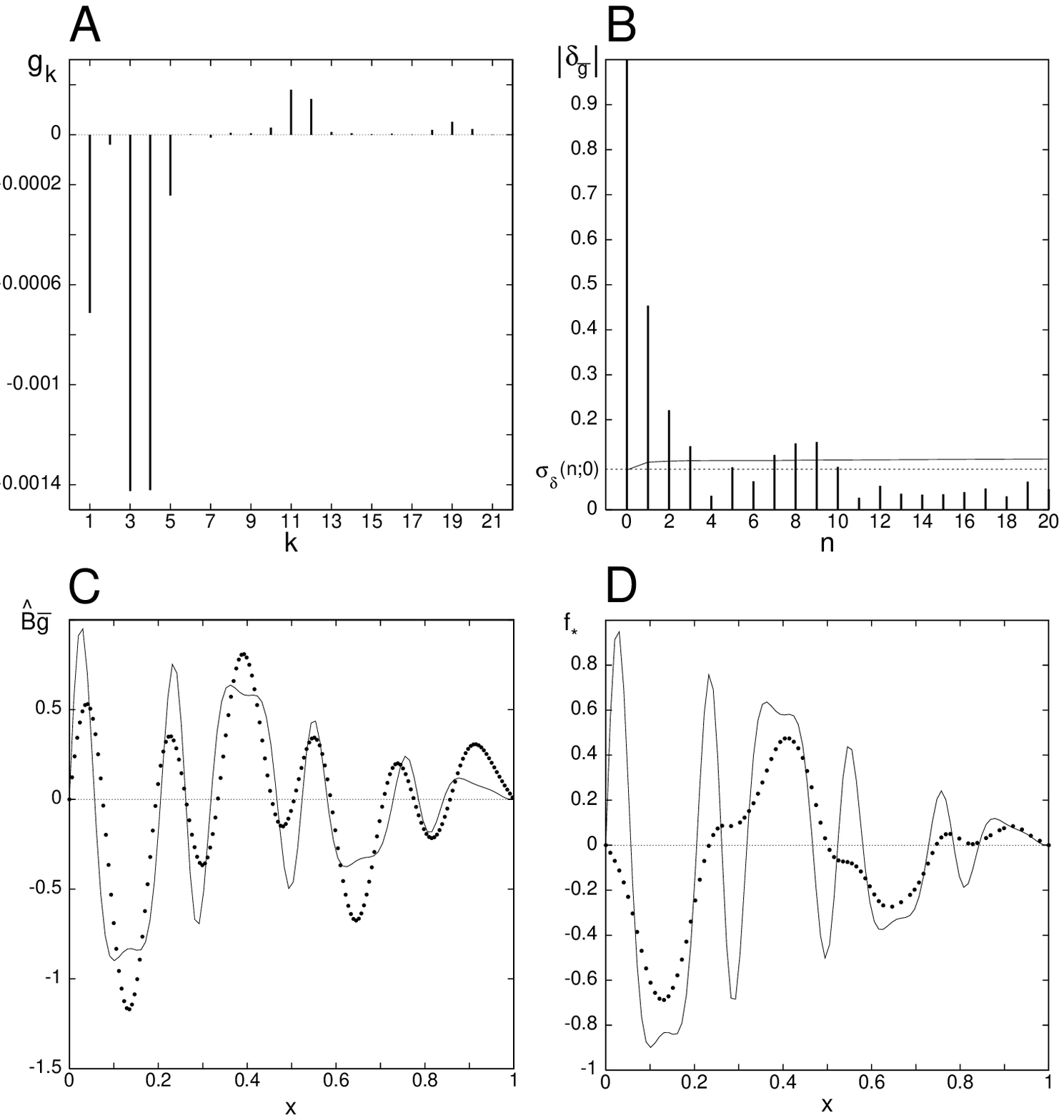}
\caption{\label{figura_3} Example 4: $f_4(x)=(1-x)\sin(5\sin(12 x))$, $\epsilon = 10^{-4}$, $\SNR \simeq 4.6 \dB$, $N = 512$.
(A) Noiseless Fourier coefficients $g_k$.
(B) Modulus of the autocorrelation function. ${\bf Q} = \{1,2,3,5,7,8,9\}$, ${\cal I}_k = \{1,3,4,9,11,12\}$.
(C) Comparison between the actual solution $f_4(x)$ (solid line) and the regularized solution
$\widehat{B}\overline{g} (x)$ (dots).
(D) Comparison between the actual solution $f_4(x)$ (solid line) and the variational solution
$f_{\star}(x)$ (see (\protect\ref{diciannove})).
}
\end{center}
\end{figure}

A more complex example is shown in Figure \ref{figura_3}. Following the trace of the previous example,
here we have the input function $f_4$ which is characterized by having the significant Fourier components
grouped in different ranges of the $k$ axis. Consequently, the Fourier coefficients $g_k$ that clearly emerge from the
noise (in this example $\epsilon=10^{-4}$) are quite sparse in the range $1\leq k \leq 12$ (see Figure \ref{figura_3}A).
The plot of the regularized solution $\widehat{B}\overline{g}$, obtained from the analysis of the
autocorrelation function shown in Figure \ref{figura_3}B, shows an acceptable agreement with
the real solution $f_4$, even though the procedure failed in detecting the coefficient at $k=5$.
On the contrary, the ``nontruncated'' (in the sense that the sum runs up to $N$) solution
$f_{\star}$ (see (\ref{diciannove})),
which is displayed in Figure \ref{figura_3}D, yields a rather poor reconstruction either because
the constraint operator $C$ smooths out too many frequencies or because distortions are
introduced by those coefficients which are essentially noise (e.g., $k=2,6,7,8,9,10$).
Of course, the variational reconstruction could be considerably improved by choosing a more appropriate operator
$C$ and different values for the parameters $\epsilon$ and $E$, but this would require
more precise a priori knowledge on the actual solution.

In conclusion, some final remarks.
The method of regularization based on the analysis of the correlation function of the data
allows to pick out the Fourier components of the noisy data which are likely to carry
exploitable information on the unknown solution, and at the same time, for rejecting the ones
dominated by the noise. Frequency selectivity is not featured by methods of regularization
that basically work as low-pass filters, and we have seen this inherent limit through examples
in which frequency selectivity is essential for a satisfactory reconstruction.

The regularized solution $\widehat{B}\overline{g}$ is founded only on a suitable analysis of the real data,
that aims at holding only the data whose information content is significant.
This approach naturally agrees with the methodology of the experimental physical science.

A moderate number of reasonable assumptions have been made in the construction of the regularized solution
$\widehat{B}\overline{g}$ (see Theorem \ref{the:2}),
and, more important, the solution itself does not depend on unknown parameters.
Even in the variational approach, methods to reduce the dependence of the solution on free parameters
have been widely investigated, and
several practical strategies for choosing the regularization parameter
$\alpha$ (see functional (\ref{nove})) have been proposed (see, for instance, \cite{Davies,Wahba}
and references therein). Since the optimal parameter is impossible to determine because the
exact solution is not known, many of these strategies can provide estimates of the asymptotically optimal
rate of convergence of the regularized solution to the real solution when the noise vanishes.

The main difficulty of the method we have proposed regards the analysis of the correlation function.
First, the correctness of the regularized solution depends on the capability of the correlation
function to catch the information content of the data and to exhibit it in an effective way.
Second, usually quite large data samples, i.e., $N$ large, are necessary in order to limit sample
fluctuations that could give rise to incorrect interpretation of the correlation function itself.


\begin{thebibliography}{20}

\bibitem{Anderson}
{\sc R. L. Anderson},
{\em Distribution of the serial correlation coefficient},
Ann. Math. Stat., 13 (1942), pp. ~1--13.

\bibitem{Balakrishnan}
{\sc A. V. Balakrishnan},
{\em Applied Functional Analysis},
Springer-Verlag, New York, 1976.

\bibitem{Bartlett}
{\sc M. S. Bartlett},
{\em Stochastic Processes: Methods and Applications},
3rd ed., Cambridge University Press, Cambridge, UK, 1978.

\bibitem{Bertero1}
{\sc M. Bertero and G. A. Viano},
{\em On probabilistic methods for the solution of improperly posed problems},
Boll. Un. Mat. Ital. B (5), 15 (1978), pp. ~483-508.

\bibitem{Bertero2}
{\sc M. Bertero, C. De Mol and G. A. Viano},
{\em On the problems of object restoration and image extrapolation in optics},
J. Math. Phys., 20 (1979), pp. ~509-521.

\bibitem{Bertero3}
{\sc M. Bertero, C. De Mol and G.A. Viano},
{\em The stability of inverse problem},
in Inverse Scattering Problems in Optics,
Springer-Verlag, Berlin, 1980, pp. ~161--212.

\bibitem{Box}
{\sc G. E. P. Box and G. M. Jenkins},
{\em Time Series Analysis},
Holden-Day, San Francisco, 1976.

\bibitem{Davies}
{\sc A. M. Davies},
{\em Optimality in regularization},
in Inverse Problems in Scattering and Imaging, M. Bertero and E.R. Pike, eds.,
Adam Hilger, Bristol, UK, 1992, pp. ~393--410.

\bibitem{Doob}
{\sc J. L. Doob},
{\em Stochastic Processes},
John Wiley, New York, 1953.

\bibitem{Engl}
{\sc H. W. Engl},
{\em Regularization methods for the stable solution of inverse problems},
Surveys Math. Indust., 3 (1993), pp. ~71-143.

\bibitem{Franklin}
{\sc J. N. Franklin},
{\em Well-posed stochastic extensions of ill-posed linear problems},
J. Math. Anal. Appl., 31 (1970), pp. ~682-716.

\bibitem{Fuller}
{\sc W. A. Fuller},
{\em Introduction to Statistical Time Series},
John Wiley, New York, 1976.

\bibitem{Gelfand}
{\sc I. M. Gel'fand and A. M. Yaglom},
{\em Calculation of tha amount of information about a random function contained in another such function},
Amer. Math. Soc. Transl. Ser. 2, 12 (1959), pp. ~199-246.

\bibitem{Groetsch}
{\sc C. W. Groetsch},
{\em The Theory of Tikhonov Regularization for Fredholm Equations of the First Kind},
Pitman, Boston, 1984.

\bibitem{Jenkins}
{\sc G. M. Jenkins and D. G. Watts},
{\em Spectral Analysis and Its Applications},
Holden-Day, San Francisco, 1968.

\bibitem{Hanke}
{\sc M. Hanke},
{\em Conjugate Gradient Type Methods for Ill-Posed Problems},
Pitman Res. Notes Math. Ser. 327,
Longman Sci. Tech., Harlow, 1995.

\bibitem{Magnoli}
{\sc N. Magnoli and G. A. Viano},
{\em On the eigenfunction expansions associated with Fredholm integral equations of first kind in presence of noise},
J. Math. Anal. Appl., 197 (1996), pp. ~188-206.

\bibitem{Maviano}
{\sc N. Magnoli and G. A. Viano},
{\em The source identification problem in electromagnetic theory},
J. Math. Phys., 38 (1997), pp. ~2366--2388.

\bibitem{Middleton}
{\sc D. Middleton},
{\em An Introduction to Statistical Communication Theory},
McGraw-Hill, New York, 1960.

\bibitem{Miller1}
{\sc K. Miller},
{\em Least square methods for ill-posed problems with a prescribed bound},
SIAM J. Math. Anal., 1 (1970), pp. ~52-74.

\bibitem{Miller2}
{\sc K. Miller and G. A. Viano},
{\em On the necessity of nearly-best-possible methods for analytic continuation of scattering data},
J. Math. Phys., 14 (1973), pp. ~1037-1047.

\bibitem{Scalas}
{\sc E. Scalas and G. A. Viano},
{\em Resolving power and information theory in signal recovery},
J. Opt. Soc. Amer. A, 10 (1993), pp. ~991-996.

\bibitem{Tikhonov}
{\sc A. Tikhonov and V. Arsenine},
{\em M\'{e}thodes de R\`{e}solution de Probl\'{e}mes Mal Pos\`{e}s},
Mir, Moscow, 1976.

\bibitem{Viano}
{\sc G. A. Viano},
{\em On the regularization of the antenna synthesis problem},
in Partial Differential Equations and Applications, P. Marcellini, G. T. Talenti, and E. Vesentini, eds.,
Marcel Dekker, 1996, pp. ~313-318.

\bibitem{Wahba}
{\sc G. Wahba},
{\em Practical approximate solutions to linear operator equations when the data are noisy},
SIAM J. Numer. Anal., 14 (1977), pp. ~651-667.

\end{thebibliography}
\end{document}